\documentclass[12pt]{article}
\usepackage{amssymb}
\usepackage{latexsym,bm}
\usepackage{graphicx}
\usepackage{amsmath}
\usepackage{mathrsfs}

\newcommand{\bea}{\begin{eqnarray*}}
\newcommand{\eea}{\end{eqnarray*}}
\newcommand{\be}{\begin{equation}}
\newcommand{\ee}{\end{equation}}
\newcommand{\ben}{\begin{eqnarray*}}
\newcommand{\een}{\end{eqnarray*}}

\date{}
\begin{document}
\title{The maximum number of $K_{r_1,\ldots,r_s}$ in graphs with a given circumference or matching number\footnote{E-mail addresses: {\tt mathdzhang@163.com}.}}
\author{\hskip -10mm Leilei Zhang\\
{\hskip -10mm \small Department of Mathematics, East China Normal University, Shanghai 200241, China}}\maketitle
\begin{abstract}
Let $K_{r_1,\ldots,r_s}$ denote the complete multipartite graph with class sizes $r_1,\ldots,r_s$ and let $K_s$ denote the complete graph of order $s$.  In 2018, Luo determined the maximum number of $K_s$ in 2-connected graphs with a given circumference. Recently, Lu, Yuan and Zhang determined the maximum number of $K_{r_1,r_2}$ in 2-connected graphs with a given circumference, and Wang determined the maximum number of $K_s$ or $K_{r_1,1,\ldots,1}$ in graphs with given matching number. Motivated by these works, we determine the maximum number of $K_{r_1,\ldots,r_s}$ in $2$-connected graphs with given circumference and large minimum degree. The maximum number of $K_{r_1,\ldots,r_s}$ in graphs with given matching number and large minimum degree is also given. Consequently, we determine the maximum number of $K_{r_1,\ldots,r_s}$ in graphs with a given circumference or matching number. We also solve the corresponding problems for
graphs with a given detour order.
\end{abstract}



\section{ Introduction}
We consider finite simple graphs, and use standard terminology and notations. Denote by $V(G)$ and $E(G)$ the vertex set and edge set of a graph $G.$ The {\it order} of a graph is its number of vertices, and the {\it size} is its number of edges. Let $e(G)$ denote the size of $G.$ For a vertex $v$ in a graph, we denote by $d(v)$ and $N(v)$ the degree of $v$ and the neighborhood of $v$ in $G,$ respectively. For $S\subseteq V(G),$ we denote by $N_S(v)$ the set $S\cap N(v)$ and denote $d_S(v)=|N_S(v)|.$ We denote by $\delta(G)$ the {\it minimum degree} of a graph $G.$  For two vertices $u$ and $v$, we use the symbol $u\leftrightarrow v$ to mean that $u$ and $v$ are adjacent and use $u\nleftrightarrow v$ to mean that $u$ and $v$ are non-adjacent. For graphs we will use equality up to isomorphism, so $G_1=G_2$ means that $G_1$ and $G_2$ are isomorphic. $\overline{G}$ denotes the complement of a graph $G.$ For two graphs $G$ and $H,$ $G\vee H$ denotes the {\it join} of $G$ and $H,$ which is obtained from the disjoint union $G+H$ by adding edges joining every vertex of $G$ to every vertex of $H.$ $G\cup H$ is the graph with vertex set $V(G)\cup V(H)$ and edge set $E(G)\cup E(H).$ We use $kG$ to denote $\bigcup_{i=1}^kG.$ For a positive integer $k,$ let $[k]:=\{1,2,\ldots,n\}$ and $k!:=1\times2\times\dots\times k.$

We use $P_k$ and $C_k$ to denote a path on $k$ vertices and a cycle on $k$ vertices, respectively. A {\it matching} $M$ is a set of pairwise nonadjacent
 edges of $G$. The {\it matching number $\alpha'(G)$} is the size of a maximum matching. $K_n$ denotes the complete graph of order $n.$ For positive integers $r_1,r_2,\ldots,r_s,$ denote by $K_{r_1,r_2,\ldots,r_s}$ the complete multipartite graph with $s$ partes of size $r_1,r_2,\ldots,r_s.$ Let $\mathcal{R}= \{r_1,r_2,\ldots,r_s\},$ and for convenience, we use $K_\mathcal{R}$ to denote $K_{r_1,r_2,\ldots,r_s}.$

A graph is called {\it Hamiltonian} if it has a Hamilton cycle; otherwise it is {\it nonhamiltonian}. In 1961 Ore \cite{8} determined the maximum size of a nonhamiltonian graph with a given order and also determined the extremal graphs.

{\bf Theorem A.} (Ore \cite{8}) {\it The maximum size of a nonhamiltonian graph of order $n$ is $\binom{n-1}{2}+1$ and this size is attained by a graph $G$ if and only if $G=K_1\vee (K_{n-2}+K_1)$ or $G=K_2\vee \overline{K_3}.$}

Bondy \cite{1} gave a new proof of Theorem A. It is natural to ask the same question by putting constraints on the graphs. In 1962 Erd\H{o}s [5] determined
the maximum size of a nonhamiltonian graph of order $n$ and minimum degree at least $k.$

{\bf Theorem B.} (Erd\H{o}s \cite{6}) {\it Let $n$, $k$ be integers with $1\le k\le \lfloor \frac{n-1}{2}\rfloor,$ and set $h(n,k)=\binom{n-k}{2}+k^2.$ If $G$ is a nonhamiltonian graph of order $n$ with $\delta(G)\ge k$, then $ e(G)\le {\rm max}\{h(n,k),\, h(n,\left\lfloor\frac{n-1}{2}\right\rfloor)\}.$ }

F\"{u}redi, Kostochka and Luo \cite{5} proved a stability version of this theorem in 2018. The {\it circumference} $c(G)$ of a graph $G$ is the length of a longest cycle in $G.$ Determining the circumference of a graph is a classical problem in graph theory. It is well known that even to determine whether a given graph is hamiltonian or not is NP-hard. One cornerstone in this direction is the following celebrated Erd\H{o}s-Gallai theorem.

{\bf Theorem C.} (Erd\H{o}s and Gallai \cite{11}) {\it Let $G$ be a graph with order $n$ and circumference $c$. Then $e(G)\le c(n-1)/2.$}

This theorem is sharp if $n-1$ is divisible by $c-1.$ This fact can be seen by considering the graph $K_1\vee (\frac{n-1}{c-1}K_{c-1}).$ Note that the order of the longest path of $G$ is $p$ if and only if the circumference of $ K_1\vee G$ is $p+1$ (\cite{10} p.292). Theorem C also implies that if an $n$-vertex graph $G$ contains no path of length $p$, then $e(G)\le (p-1)n/2.$  F\"{u}redi, Kostochka and Verstra\"{e}te [9] proved a stability version of Theorem C.

{\bf Notation 1.} Let $r=\sum_{i=1}^sr_i.$ Denote by $P(r;r_1,r_2,\ldots,r_s)$ the multinomial coefficients, i.e. $P(r;r_1,r_2,\ldots,r_s)=\binom{r}{r_1\ r_2\ \cdots\ r_s}=\frac{r!}{r_1!r_2!\cdots r_s!}$

{\bf Notation 2.} Fix $n-1 \geq c\geq 2k.$ Let $G(n,c,k)=K_k\vee(K_{c+1-2k}+\overline{K_{n-c-1+k}})$ and $r=\sum_{i=1}^sr_i$. Suppose $\mathcal{R}= \{r_1,r_2,\ldots,r_s\}$ contains $q$ pairwise different $r'_i$ and let $m_i$ be the number of times $r'_i$ appears in set $\mathcal{R}.$ Denote by $g_\mathcal{R}(n,c,k)$ the number of $K_{\mathcal{R}}$ in $G(n,c,k)$; more precisely,
\begin{small}
\begin{align*}
g_\mathcal{R}(n,c,k)
=&\sum_{i=1}^{n-c-1+k}\sum_{j=1}^{s}\frac{P(r-r_j;\mathcal{R}\backslash\{r_j\})}{\Pi_{l=1}^{q}m_l!}\binom{k}{r-r_j}\binom{n-r+r_j-i}{r_j-1}
 +\frac{P(r;\mathcal{R})}{\Pi_{l=1}^{q}m_l!}\binom{c+1-k}{r}.
\end{align*}
\end{small}

For $G(n,c,k)=K_k\vee(K_{c+1-2k}+\overline{K_{n-c-1+k}}),$ remark the vertices in $\overline{K_{n-c-1+k}}$ by $x_1,\dots,x_{n-c-1+k}.$ Note that for $i\in [n-c-1+k],$ the number of copies of $K_\mathcal{R}$ containing $x_i$ in $G(n,c,k)-\{x_1,\ldots,x_{i-1}\}$ is
$\sum_{j=1}^{s}\frac{P(r-r_j;\mathcal{R}\backslash\{r_j\})}{\Pi_{l=1}^{q}m_l!}\binom{k}{r-r_j}\binom{n-r+r_j-i}{r_j-1}.$ The number of copies of $K_\mathcal{R}$ in $K_k\vee K_{c+1-2k}$ is $\frac{P(r;\mathcal{R})}{\Pi_{l=1}^{q}m_l!}\binom{c+1-k}{r}.$ Hence the number of copies of $K_\mathcal{R}$ in $G(n,c,k)$ is $g_\mathcal{R}(n,c,k).$ For $r_1=\ldots=r_s=1,$ we write $g_s(n,c,k)$ for $g_\mathcal{R}(n,c,k)$ which equals to the number of $K_s$ in $G(n,c,k).$

{\bf Conjecture D.} (Woodall \cite{13}) {\it Let $n-1\ge c\ge 4$. The maximum size of a $2$-connected graph with order $n,$ circumference $c$ and $\delta(G)\ge k$ is
${\rm max}\{g_2(n,c,k), \, g_2(n,c,\lfloor c/2 \rfloor)\}.$}

Woodall's conjecture D has been proved (see \cite{14,9,7,12}). For further developments on this topic, see \cite{23,3,4,22,16,2}. Improving Theorem C, Kopylov \cite{9} proved the following result in 1977.

{\bf Theorem E.} (Kopylov \cite{9}) {\it The maximum size of a $2$-connected nonhamintonian graph of order $n$ and circumference $c$ is
${\rm max}\{g_2(n,c,2),\, g_2(n,c,\lfloor c/2\rfloor)\}.$}

{\bf Definition 1.} The order of a longest path in a graph $G$ is called the detour order of $G.$

A graph is called {\it traceable} if it has a Hamilton path; otherwise it is {\it nontraceable}. For the detour, Kopylov \cite{9} and Balister, Gyori, Lehel, Schelp \cite {17} independently gave the maximum size of a connected graph with given detour order.

{\bf Notation 3.} For complete multipartite graph $K_{r_1,\ldots,r_s}$, we use $N(K_{r_1,\ldots,r_s},G)$ to denote the number of $K_{r_1,\ldots,r_s}$ in $G$. If $K_{r_1,\ldots,r_s}$ is a complete graph, for convenience, we denote by $N(K_{s},G)$ the number of $K_{s}$ in $G.$

Generalizing Theorem E (but using Kopylov's proof idea), Luo \cite{7} proved the following result.

{\bf Theorem F} (Luo \cite{7}). {\it Let $n-1\geq c\geq 4$ and let $s\ge 2$. If $G$ is a $2$-connected $n$-vertex graph with $c(G)=c$, then
$ N(K_s,G)\leq {\rm max}\{g_s(n,c,2),\, g_s(n,c,\lfloor c/2\rfloor)\}.$}

Ning and Peng presented an extension of Theorem F by adding a restraint on the minimum degree.

{\bf Theorem G} (Ning and Peng \cite{12}). {\it Let $n-1\ge c$ and $s\geq 2$. If $G$ is a $2$-connected graph of order $n$ with  $c(G)=c$ and minimum degree $\delta(G)\geq k$, then $ N(K_s,G)\leq {\rm max}\{g_s(n,c,k),\, g_s(n,c,\lfloor c/2\rfloor)\}.$}

Denote by $S_r$ the star on $r$ vertices and $C_r$ the cycle on $r$ vertices, respectively. Gy\H{o}ri, Salia, Tompkins, Zamora \cite{18} considered the maximum number of $S_r$ and $C_4$ copies in $n$-vertex graphs with given circumference or detour order for sufficiently larger $n$.

{\bf Theorem H} (Gy\H{o}ri, Salia, Tompkins, Zamora \cite{18}). {\it Let $n-1\ge c\ge 4,$ $r\ge3$ and for sufficiently large $n$. If $G$ is a $n$-vertex graph with circumference $c,$ then $ N(S_r,G)\le g_{1,r-1}(n,c,\lfloor c/2\rfloor).$}

{\bf Theorem I} (Gy\H{o}ri, Salia, Tompkins, Zamora \cite{18}). {\it Let $n-1\ge c\ge 4,$ $r\ge3$ and for sufficiently large $n$. If $G$ is a $n$-vertex graph with circumference $c,$ then $ N(C_4,G)\le g_{2,2}(n,c,\lfloor c/2\rfloor).$}

Note that $S_r$ and $C_4$ are special complete bipartite graphs. Lu, Yuan and Zhang \cite{19} generalized the above results to all complete bipartite graphs.

{\bf Theorem J} (Lu, Yuan, Zhang \cite{19}). {\it Let $n-1\geq c\geq 4$. If $G$ is a $2$-connected $n$-vertex graph with circumference $c$, then
$ N(K_{r_1,r_2},G)\leq {\rm max}\{g_{r_1,r_2}(n,c,2),\, g_{r_1,r_2}(n,c,\lfloor c/2\rfloor)\}.$}

{\bf Theorem K} (Lu, Yuan, Zhang \cite{19}). {\it Let $n-1\geq p\geq 3$. If $G$ is a connected $n$-vertex graph with detour order $p$, then
$ N(K_{r_1,r_2},G)\leq {\rm max}\{g_{r_1,r_2}(n,p-1,1),\, g_{r_1,r_2}(n,p-1,\lfloor (p-1)/2\rfloor)\}.$}

In the same paper \cite{19}, for $r_1\ge r_2 \ge 2,$ they also determined the maximum number of copies of $K_{r_1,r_2}$ in $n$-vertex graphs with given detour order. Recently, Duan, Ni, Peng, Wang, and Yang \cite{21} determined the maximum number of $K_s$ in $n$-vertex graph $G$ with given matching number $\alpha'(G)$ and $\delta (G)\ge k.$
Let $K_{r,s}^*$ be obtained by taking a copy $K_{r,s}$ and joining each pair of vertices inside the part with size $r.$ By using the shifting method, Wang \cite{20} gave the following results.

{\bf Theorem L} (Wang \cite{20}). {\it Let $s\ge 2$ and $n\ge 2k+1$. If $G$ is a $n$-vertex graph with matching number $k,$ then $ N(K_{s},G)\leq {\rm max}\{\binom{2k+1}{s},\binom{k}{s}+(n-k)\binom{k}{s-1}\}$}

{\bf Theorem M} (Wang \cite{20}). {\it Let $r\ge 1$, $s\ge 2$ and $n\ge 2k+1$. If $G$ is a $n$-vertex graph with matching number $k,$ then $ N(K^*_{r,s},G)\leq {\rm max}\{\binom{2k+1}{s+r}\binom{s+r}{s},\binom{k}{r}\binom{n-s}{s}+(n-k)\binom{k}{s+r-1}\binom{s+r-1}{s}\}$  }

Note that $K_s$, $K_{r_1,r_2}$ and $K_{r,s}^*$ are special complete multipartite graphs.  Motivated by these works, we present an extension of Theorem H$\sim$M by adding a restraint on the minimum degree. We determine the maximum number of $K_{r_1,\ldots,r_s}$ in $2$-connected graphs with given circumference and large minimum degree. This result generalizes Theorem H. The maximum number of $K_{r_1,\ldots,r_s}$ in graphs with given matching number and large minimum degree is also given. This result generalizes Theorem L and Theorem M. The maximum number of copies of $K_{r_1,\ldots,r_s}$ in $n$-vertex connected graph with given detour order and large minimum degree is also determined. This result generalizes Theorem K.  Moreover, Lu, Yuan and Zhang \cite{19} determined the maximum number of copies of $K_{r_1,r_2}$ in $n$-vertex graphs with a given detour order for $r_1\ge r_2\ge 2$. As corollaries of our main results, we also determine the maximum number of $K_{r_1,\ldots,r_s}$ in $n$-vertex graph with a given circumference or detour order or matching number for all $r_i\ge 1, \ i\in[s]$. This result generalizes Theorem H and Theorem I.

{\bf Theorem 1.} {\it Let $G$ be a 2-connected nonhamiltonian graph of order $n,$ circumference $c$ and minimum degree at least $k.$ Then
$$N(K_{r_1,\ldots,r_s},G)\leq {\rm max}\{g_{r_1,\ldots,r_s}(n,c,k),\, g_{r_1,\ldots,r_s}(n,c,\lfloor c/2\rfloor)\}.$$}
This theorem is sharp as shown by the examples $G(n,c,k)$ and $G(n,c,\lfloor c/2\rfloor).$

{\bf Theorem 2.} {\it Let $G$ be a connected nontraceable graph of order $n,$ detour order $p$ and minimum degree at least $k.$ Then
$$N(K_{r_1,\ldots,r_s},G)\leq {\rm max}\{g_{r_1,\ldots,r_s}(n,p-1,k),\, g_{r_1,\ldots,r_s}(n,p-1,\lfloor (p-1)/2\rfloor)\}.$$}
Again, this theorem is sharp as shown by the examples $G(n,p-1,k)$ and $G(n,p-1,\lfloor (p-1)/2\rfloor))$

{\bf Theorem 3.} {\it Let $G$ be a graph of of order $n\geq 2\alpha'+2,$ matching number $\alpha'$ and minimum degree at least $k.$ Then
$$N(K_{r_1,\ldots,r_s},G)\leq {\rm max}\{g_{r_1,\ldots,r_s}(n,2\alpha',k),\, g_{r_1,\ldots,r_s}(n,2\alpha',\alpha')\}.$$}

The rest of the paper is organized as follows. In Section 2, we will prove Theorem 1, Theorem 2 and Theorem 3. In Section 3 we give some corollaries. Concluding remarks will be given in Section 4.

\section{Proof of the main results}
To prove Theorem 1, we will need the following lemmas and definitions.

{\bf Lemma 4.} {\it $g(k)$ is a convex function of $k$ if and only if $g(k+1)+g(k-1)-2g(k)\ge 0$.}

{\bf Lemma 5.} {\it $g_{r_1,\dots,r_s}(n,c,k)$ is a convex function of $k$.}

{\bf Proof.} Suppose $\mathcal{R}= \{r_1,r_2,\ldots,r_s\}$ contains $q$ pairwise different $r'_i$ and let $m_i$ be the number of times $r'_i$ appears in set $\mathcal{R}.$ Denote by $r=\sum_{i=1}^sr_i=\sum_{i=1}^qm_ir_i'.$ Remained that $g_{\mathcal{R}}(n,c,k)$ is the number of $K_\mathcal{R}$ in $G(n,c,k)=K_k\vee(K_{c+1-2k}+\overline{K_{n-c-1+k}}).$ Note that the number of copies of $K_{\mathcal{R}}$ inside $K_k\vee K_{c+1-2k}$ are $\frac{P(r;\mathcal{R})}{\Pi_{l=1}^{q}m_l!}\binom{c+1-k}{r}.$ The number of copies of $K_{\mathcal{R}}$ not inside $K_k\vee K_{c+1-2k}$ are $\sum_{i=1}^s\frac{P(r-r_i;\mathcal{R}\setminus r_i)}{\Pi_{l=1}^{q}m_l!}\binom{k}{r-r_i}\left[\binom{n-r+r_i}{r_i}-\binom{c+1-k-r+r_i}{r_i}\right]$. Then we have
\begin{align*}
g_{\mathcal{R}}(n,c,k)=&\frac{P(r;\mathcal{R})}{\Pi_{l=1}^{q}m_l!}\binom{c+1-k}{r}+\sum_{i=1}^s\frac{P(r-r_i;\mathcal{R}\setminus r_i)}{\Pi_{l=1}^{q}m_l!}\binom{k}{r-r_i}\binom{n-r+r_i}{r_i}\\
&-\sum_{i=1}^s\frac{P(r-r_i;\mathcal{R}\setminus r_i)}{\Pi_{l=1}^{q}m_l!}\binom{k}{r-r_i}\binom{c+1-k-r+r_i}{r_i}.
\end{align*}
Not that
\begin{align*}
&g_{\mathcal{R}}(n,c,k+1)-g_{\mathcal{R}}(n,c,k)\\
=&-\frac{P(r;\mathcal{R})}{\Pi_{l=1}^{q}m_l!}\binom{c-k}{r-1}+\sum_{i=1}^s\frac{P(r-r_i;\mathcal{R}\setminus r_i)}{\Pi_{l=1}^{q}m_l!}\binom{k}{r-r_i-1}\binom{n-r+r_i}{r_i}\\
&+\sum_{i=1}^s\frac{P(r-r_i;\mathcal{R}\setminus r_i)}{\Pi_{l=1}^{q}m_l!}\left[\binom{k}{r-r_i}\binom{c-k-r+r_i}{r_i-1}-\binom{k}{r-r_i-1}\binom{c-k-r+r_i}{r_i}\right].
\end{align*}
We have
\begin{align*}
& g_{\mathcal{R}}(n,c,k+1)+g_{\mathcal{R}}(n,c,k-1)-2g_{\mathcal{R}}(n,c,k)\\
\ge&\frac{P(r;\mathcal{R})}{\Pi_{l=1}^{q}m_l!}\binom{c-k}{r-2}+\sum_{i=1}^s\frac{P(r-r_i;\mathcal{R}\setminus r_i)}{\Pi_{l=1}^{q}m_l!}\binom{k-1}{r-r_i-2}\binom{n-r+r_i}{r_i}\\
&-\sum_{i=1}^s\frac{P(r-r_i;\mathcal{R}\setminus r_i)}{\Pi_{l=1}^{q}m_l!}\left[\binom{k-1}{r-r_i-2}\binom{c-k-r+r_i}{r_i}+\binom{k-1}{r-r_i}\binom{c-k-r+r_i}{r_i-2}\right]\\
\ge& \frac{P(r;\mathcal{R})}{\Pi_{l=1}^{q}m_l!}\binom{c-k}{r-2}-\sum_{i=1}^s\frac{P(r-r_i;\mathcal{R}\setminus r_i)}{\Pi_{l=1}^{q}m_l!}\binom{k-1}{r-r_i}\binom{c-k-r+r_i}{r_i-1}.
\end{align*}
By lemma 4, it is enough to prove the following inequality.
\begin{align*} \frac{P(r;\mathcal{R})}{\Pi_{l=1}^{q}m_l!}\binom{c-k}{r-2}\ge\sum_{i=1}^s\frac{P(r-r_i;\mathcal{R}\setminus r_i)}{\Pi_{l=1}^{q}m_l!}\binom{k-1}{r-r_i}\binom{c-k-r+r_i}{r_i-1}.
\end{align*}
By simplification, we have
\begin{align*}
1\ge\sum_{i=1}^s\frac{r_i(r_i-1)}{r(r-1)}\frac{(k-1)(k-2)\ldots(k-r+r_i)}{(c-k)(c-k-1)\ldots(c-k-r+r_i+1)}
\end{align*}
Since $c\ge 2k$ and $r=\sum_{i=1}^sr_i,$ we have
\begin{align*}
\sum_{i=1}^s\frac{r_i(r_i-1)}{r(r-1)}\frac{(k-1)(k-2)\ldots(k-r+r_i)}{(c-k)(c-k-1)\ldots(c-k-r+r_i+1)}<\sum_{i=1}^s \frac{r_i(r_i-1)}{r(r-1)}\le 1
\end{align*}
 This completes the proof.
 
We remark that there have been a lot of works exploiting the convexity property of an appropriate graph-related function. However in most of these works convexity is trivial to check. We present in this paper for the first time a detailed proof of the important observation that $g_{r_1,\dots,r_s}(n,c,k)$ (see Notation 2) is convex.

{\bf Lemma 6.} (Kopylov \cite{9}) {\it Let $P$ be an $(x,y)$-path of length $m$ in a $2$-connected graph $G$. Then
$c(G)\geq {\rm min}\{m+1,\, d_P(x)+d_P(y)\}.$}

The following definition of {\it$t$-disintegration} of a graph is due to Kopylov \cite{9}.

{\bf Definition 2.} ($t$-disintegration of a graph, Kopylov \cite{9}). Let $G$ be a graph and $t$ be a natural number. Delete all vertices of degree at most $t$ from $G$; for the resulting graph $G'$, we again delete all vertices of degree at most $t$ from $G'$. Iterating this process until we finally obtain a graph, denoted by $D(G; t)$, such that either $D(G; t)$ is a null graph or $\delta(D(G; t))\geq t+1.$ The graph $D(G; t)$ is called the {\it $(t+1)$-core of $G.$}

{\bf Definition 3.} $G$ is called {\it edge-maximal} with respect to the circumference if for any $e\in E(\overline{G}),$ $c(G+e)> c(G).$

Now we are ready to prove Theorem 8.  We will use ideas from \cite{9} (proof of Theorem 3), \cite{7} (proof of Theorem 1.4) and \cite{12} (proof of Theorem 3.4).
We also need to treat new situations, since more precise conditions are given in our problem.

{\bf Proof of Theorem 1.} Let $n-1\geq c\geq 2k ,t=\lfloor c/2\rfloor.$ It is easy to verify that the graphs $G(n,c,k)$ and $G(n,c,t)$ stated in Notation 1 are graphs of order $n$, circumference $c$ and minimum degree at least $k.$ The number of $K_\mathcal{R}$ in $G(n,c,k)$ or $G(n,c,t)$ is $g_{\mathcal{R}}(n,c,k)$ or $g_{\mathcal{R}}(n,c,t).$

Suppose $Q$ is edge-maximal with respect to the circumference $c.$ Thus each pair of non-adjacent vertices in $Q$ is connected by a path of length at least $c.$ Let $D(Q;t)$ denote the $(t+1)$-core of $Q$, i.e., the resulting graph of applying $t$-disintegration to $Q.$ For convenience, let $D=D(Q;t)$. We distinguish two cases.

{\bf Case 1.} $D$ is a null graph. In the $t$-disintegration process, denote $Q_1=Q$ and $Q_{i+1}=Q_i-x_i,1\leq i\leq n-1$ where $x_i$ is a vertex of degree at most $t$ in $Q_i$. For the first $n-t$ vertices, once $x_i$ is deleted, we delete at most $\sum_{j=1}^{s}\frac{P(r;\mathcal{R}\backslash \{r_j\})}{\Pi_{l=1}^{q}m_l!}
\binom{t}{r}\binom{n-r+r_j-i}{r_j-1}$ copies of $K_\mathcal{R}.$ The number of $K_\mathcal{R}$ contained in the last $t$ vertices is at most $\frac{P(r;\mathcal{R})}{\Pi_{l=1}^{q}m_l!}\binom{t}{r}.$ Thus
\begin{align*}
N(K_\mathcal{R},Q)
&\leq \sum_{i=1}^{n-t} \sum_{j=1}^{s}\frac{P(r;\mathcal{R}\backslash\{r_j\})}{\Pi_{l=1}^{q}m_l!}\binom{t}{r}\binom{n-r+r_j-i}{r_j-1}
 +\frac{P(r;\mathcal{R})}{\Pi_{l=1}^{q}m_l!}\binom{t}{r}\\
&\leq g_\mathcal{R}(n,c,t).
\end{align*}

{\bf Case 2.} $D$ is not a null graph. Let $d=|D|$, we claim that $D$ is a complete graph and $\delta (Q)=k \leq c-d+1$.

If there exist two vertices that are not adjacent in $D$, then in $Q$, there is a path of length at least $c$ with these vertices as its endpoints. Among all nonadjacent pairs of vertices in $D$, choose $x, y\in V(D)$ such that $|P(x,y)|={\rm max}\{|P(u,v)|:u,v\in V(D), u\nleftrightarrow v\}.$ Let $P_1=P(x,y)$. We next show $N_D(x)\subseteq V(P_1)$ and $N_D(y)\subseteq V(P_1).$ If $x$ has a neighbor $z\in V(D)$ and $z\notin V(P_1)$ , then either $yz\in E(Q)$ and $zx\cup P_1\cup yz$ is a cycle of length at least $c+1$, or $yz\notin E(Q)$ and so $P_1\cup xz$ is a longer path. This contradicts the maximality of $P_1.$ Similarly for $y$, we have $N_D(y)\subseteq V(P_1).$ Hence, by Lemma 6, $Q$ has a cycle of length at least $min\{c+1,d_{P_1}(x)+d_{P_1}(y)\}\geq min\{c+1,2(t+1)\}= c+1,$ contradiction. Thus $D$ is a complete graph.

Suppose $k\geq c-d+2$. Then $d\geq c-k+2.$ By the definition of $t$-disintegration, the minimum degree of $D$ is at least $t+1$, so we have $d\geq t+2.$ If $u\in V(Q)\backslash V(D)$, then $u$ is not adjacent to at least one vertex in $D.$  Choose $x\in V(Q)\backslash V(D)$ and $y\in V(D)$ such that $|P(x,y)|={\rm max}\{|P(u,v)|:u\in V(Q)\backslash V(D), v\in V(D), u\nleftrightarrow v\}.$ Denote $P_2=P(x,y)$. Since each pair of non-adjacent vertices in $Q$ is connected by a path of length at least $c,$ by the maximality of $P_2,$ we have the size of $P_2$ is at least $c.$ We claim that $N_Q(x)\subseteq V(P_2)$ and $N_D(y)\subseteq V(P_2).$  If $x$ has a neighbor $z\in V(D)$ and $z\notin V(P_2).$ Since $D$ is a complete graph, $zx\cup P_2\cup yz$ is a cycle of length at least $c+1$, contradiction. If $x$ has a neighbor $z\in V(Q)\backslash V(D)$ and $z\notin V(P_2),$ then either $yz\in E(Q)$ and $zx\cup P_2\cup yz$ is a cycle of length at least $c+1$, or $yz\notin E(Q)$ and so $P_2\cup xz$ is a longer path. This contradicts the maximality of $P_2.$ Similarly for $y$, we have $N_D(y)\subseteq V(P_2).$ Hence, by Lemma 6, there is a cycle with length at least ${\rm min}\{c+1,d_{P_2}(x)+d_{P_2}(y)\}\geq {\rm min}\{c+1,k+d-1\}\geq {\rm min}\{c+1,k+c-k+1\}= c+1,$ where the second inequality follows from $d\geq c-k+2,$  contradiction. Thus $k\leq c-d+1.$

Apply $(c-d+1)$-disintegration to $Q$, and let $D'=D(Q;c-d+1)$ be the resulting graph. Recall that $d\geq t+2$. We have $k\leq c-d+1\leq c-t-1 \leq t$. Then $D\subseteq D'.$ There are two cases.

(a) If $D'=D$, then $|D'|=|D|=d.$ By the definition of $(c-d+1)$-disintegration, we have
\begin{align*}
N(K_\mathcal{R},Q)
&\leq \sum_{i=1}^{n-d} \sum_{j=1}^{s}\frac{P(r-r_j;\mathcal{R}\backslash \{r_j\})}{\Pi_{l=1}^{q}m_l!}\binom{c-d+1}{r-r_j}\binom{n-r+r_j-i}{r_j-1}
  +\frac{P(r;\mathcal{R})}{\Pi_{l=1}^{q}m_l!}\binom{d}{r}\\
&= g_\mathcal{R}(n,c,c-d+1)\\
&\leq {\rm max}\{ g_{\mathcal{R}}(n,c,k),g_{\mathcal{R}}(n,c,t)\}.
\end{align*}
The third inequality follows from the condition $k\leq c-d+1 \leq t$ and that the function $g_\mathcal{R}(n,c,x)$ is convex for $x\in[k,t].$

(b) Otherwise $D\neq D'.$ If $u\in V(D')\backslash V(D)$, then $u$ is not adjacent to at least one vertex in $D.$  Among all these nonadjacent pairs of vertices, choose $x\in V(D')\backslash V(D)$ and $y\in V(D)$ such that $|P(x,y)|={\rm max}\{|P(u,v)|: u\in V(D')\backslash V(D), v\in V(D), u\nleftrightarrow v\}.$ Denote $P_3=P(x,y)$ for convenience. As before, we have  $N_{D'}(x)\subseteq V(P_3)$ and $N_D(y)\subseteq V(P_3).$  By Lemma 6, there is a cycle with length at least ${\rm min}\{c+1, d_{P_3}(x)+d_{P_3}(y)\} \geq {\rm min}\{c+1,(c-d+2)+(d-1)\}=c+1$, where the second inequality follows from the condition $D'$ is the $(c-d+2)$-core of $Q$ and $|D|=d,$ contradiction.

This completes the proof. \hfill $\Box$

The proof of Theorem 2 follows the same steps as the proof of theorem 1. So we will omit some details.

{\bf Proof of Theorem 2.} Let $n-1\geq p\geq 2k+1 ,t=\lfloor (p-1)/2\rfloor.$ It is easy to verify that the graphs $G(n,p-1,k)$ and $G(n,p-1,t)$ stated in Notation 1 are graphs of order $n$, detour order $p$ and minimum degree at least $k.$ The number of $K_\mathcal{R}$ in $G(n,p-1,k)$ or $G(n,p-1,t)$ is $g_{\mathcal{R}}(n,p-1,k)$ or $g_{\mathcal{R}}(n,p-1,t).$

Suppose $G$ is a $n$-vertex graph detour order $p$ and $\delta(G)\ge k.$ Let $Q=G\vee \{v_0\}.$ Add edges to $Q$ until addition any edge creates a cycle of length at least $p+1.$ Denote the resulting graph by $Q'.$ Thus each pair of non-adjacent vertices in $Q'$ is connected by a path of length at least $p+2.$ Apply $(t+1)$-disintegration to $Q',$ where if necessary, we delete $v_0$ last. Let $D(Q';t+1)$ denote the $(t+2)$-core of $Q'$. We distinguish two cases.

{\bf Case 1.} $D(Q';t+1)$ is empty. For the first $n-t$ vertices in the process of $t+1$-disintegration, each of them has at most $t$ neighbors that are not $v_0.$ As in the proof of Theorem 1. Once $x_i$ is deleted, we delete at most  $\sum_{j=1}^{s}\frac{P(r-r_j;\mathcal{R}\backslash\{r_j\})}{\Pi_{l=1}^{q}m_l!}
\binom{t}{r-r_j}\binom{n-r+r_j-i}{r_j-1}$ copies of $K_\mathcal{R}$ in $Q$. Thus
\begin{align*}
N(K_\mathcal{R},Q)
&\leq \sum_{i=1}^{n-t}\sum_{j=1}^{s}\frac{P(r-r_j;\mathcal{R}\backslash\{r_j\})}{\Pi_{l=1}^{q}m_l!}\binom{t}{r-r_j}\binom{n-r+r_j-i}{r_j-1}
 +\frac{P(r;\mathcal{R})}{\Pi_{l=1}^{q}m_l!}\binom{t}{r}\\
&\leq g_\mathcal{R}(n,p-1,t).
\end{align*}

{\bf Case 2.} $D(Q';t+1)$ is not empty. The same argument as in the proof of theorem 1 also shows that $D(Q';t+1)$ is a clique, otherwise there would be a cycle of length at least $2(t+2)\geq p+2$ in $Q'.$ Note that $v_0$ must be contained in $D(Q';t+1)$ as it is  adjacent to all other vertices of $Q'.$ Let $d=|D(Q';t+1)|.$ Since very vertex in  $D(Q';t+1)$ has degree at least $t+2,$ $d \geq t+3$. By the same argument as in the proof of theorem 1,  we have $k+1\leq p-d+2\leq p-t-1\leq t+1.$

Apply $(p-d+2)$-disintegration to $Q'.$ If $D(Q';t+1)\neq D(Q';p-d+2),$ we will find a cycle of length at least $p+2$ in $Q'.$ Otherwise $D(Q';t+1)= D(Q';p-d+2).$ In $D(Q';t+1)$ the number of $K_\mathcal{R}$ that not include $v_0$ is $\frac{P(r;\mathcal{R})}{\Pi_{l=1}^{q}m_l!}\binom{d-1}{r}.$ In $Q'-D(Q';p-d+2),$ by the definition of $(p-d+2)$-disintegration, every vertex had at most $p-d+1$ neighbors that were not $v_0.$ Then
\begin{align*}
N(K_\mathcal{R},Q)
&\leq \sum_{i=1}^{n+1-d} \sum_{j=1}^{s}\frac{P(r-r_j;\mathcal{R}\backslash \{r_j\})}{\Pi_{l=1}^{q}m_l!}\binom{p-d+1}{r-r_j}\binom{n-r+r_j-i}{r_j-1}
  +\frac{P(r;\mathcal{R})}{\Pi_{l=1}^{q}m_l!}\binom{d-1}{r}\\
&= g_\mathcal{R}(n,p-1,p-d+1)\\
&\leq {\rm max}\{ g_{\mathcal{R}}(n,p-1,k),g_{\mathcal{R}}(n,p-1,t)\}.
\end{align*}
The third inequality follows from the condition $k\leq p-d+1 \leq t$ and that the function $g_\mathcal{R}(n,p-1,x)$ is convex for $x\in[k,t].$ This completes the proof. \hfill $\Box$

{\bf Definition 4.} The {\it $t$-closure} of $G$ is the graph obtained from $G$ by iteratively joining nonadjacent vertices with degree sum at least $t$ until there is no more such a pair of vertices.

{\bf Lemma 7.} (Bondy and Chv\'{a}tal \cite{9}) Let $G$ be a graph and $G'$ be the $(2k+1)$-closure of $G.$ Then $\alpha'(G')\geq k+1$ implies $\alpha'(G)\geq k+1.$

{\bf Definition 5.} $G$ is called {\it edge-maximal} with respect to the matching number if for any $e\in E(\overline{G}),$ $\alpha'(G+e)> \alpha'(G).$

{\bf Proof of Theorem 3.} It is easy to verify that the graphs $G(n,2\alpha',k)$ and $G(n,2\alpha',\alpha')$ stated in Notation 1 are graphs of order $n$, matching number $\alpha'$ and minimum degree at least $k.$ The number of $K_\mathcal{R}$ in $G(n,2\alpha',k)$ or $G(n,2\alpha',\alpha')$ is $g_{\mathcal{R}}(n,2\alpha',k)$ or $g_{\mathcal{R}}(n,2\alpha',\alpha').$

Let $G$ be a $n$-vertex graph with $\alpha'(G)=\alpha'$ and $\delta(G)\ge k.$  Let $Q$ be the $(2\alpha'+1)$-closure of $G.$ If $\alpha'(Q)\ge \alpha'+1,$ then $\alpha'(G)\ge \alpha'+1$ by lemma 7. Contradiction. Thus $\alpha'(Q)=\alpha$ and $\delta(Q)\ge \delta(G)\ge k.$ Let $D=D(Q;\alpha')$ denote the $(\alpha'+1)$-core of $Q.$ We distinguish two cases.

{\bf Case 1.} $D$ is a null graph. In the $\alpha'$-disintegration process, denote $Q_1=Q$ and $Q_{i+1}=Q_i-x_i,1\leq i\leq n-1$ where $x_i$ is a vertex of degree at most $\alpha'$ in $Q_i$. We have $d_{Q_i}(x_i)\leq \alpha',0\leq i\leq n-\alpha'-1.$ For the first $n-\alpha'-1$ vertices, once $x_i$ is deleted, we deleted at most $\sum_{j=1}^{s}\frac{P(r-r_j;\mathcal{\mathcal{R}}\backslash\{r_j\})}{\Pi_{l=1}^{q}m_l!}\binom{\alpha'}{r-r_j}\binom{n-r+r_j-i}{r_j-1}$ copies of $K_\mathcal{R}.$ The number of $K_\mathcal{R}$ contained in the least $\alpha'+1$ vertices is at most $\frac{P(r;\mathcal{R})}{\Pi_{l=1}^{\alpha'}m_l!}\binom{\alpha'+1}{r}.$ Thus
\begin{align*}
N(K_\mathcal{R},Q)
&\leq \sum_{i=1}^{n-\alpha'-1} \sum_{j=1}^{s}\frac{P(r-r_j;\mathcal{R}\backslash\{r_j\})}{\Pi_{l=1}^{q}m_l!}\binom{\alpha'}{r-r_j}\binom{n-r+r_j-i}{r_j-1}
 +\frac{P(r;\mathcal{R})}{\Pi_{l=1}^{q}m_l!}\binom{\alpha'+1}{r}\\
&=g_\mathcal{R}(n,2\alpha',\alpha').
\end{align*}

{\bf Case 2.} $D$ is not a null graph. Let $d=|D|$, we claim that $D$ is a complete graph and $k\le \delta (Q) \le 2\alpha'-d+1$.

For all $u, v\in D$, we note $d_Q(u)\geq \alpha'+1$, $d_Q(v)\geq \alpha'+1$. Since every nonadjacent vertices with degree sum at most $2\alpha'$ in $Q$ and $d_Q(u)+d_Q(v)\geq 2\alpha'+2$, $u$ and $v$ are adjacent in $D$, i.e., $D$ is a clique.

We next show $d\leq 2\alpha'-k+1.$ As $D$ is a clique and $d_D(u)\geq \alpha'+1$ for all $u\in V(D),$ we get $d\geq \alpha'+2.$ Thus, every vertex in $ V(Q)\backslash V(D)$ is not adjacent to at least two vertices in $D.$ Suppose $d\geq 2\alpha'-k+2.$ Then $d_Q(u)\geq d_D(u)\geq 2\alpha'-k+1$ for all $u\in V(D)$. Let $v$ be a vertex from $V(Q)\backslash V(D)$ which is not adjacent to $w\in V(D)$. Note that $d_Q(v)\geq k$ and $d_Q(w)+d_Q(v)\geq 2\alpha'-k+1+k=2\alpha'+1,$ which is a contradiction to the choice of $u$ and $w$. Thus $\alpha'+ 2 \leq d \leq 2\alpha'-k+1$.

Let $D'$ be the $(2\alpha'-d+2)$-core of $Q,$ i.e., the resulting graph of applying $2\alpha'-d+1$-disintegration to $Q.$ Since $d\geq \alpha'+2$, we obtain $2\alpha'-d+1\leq \alpha'.$ Therefore, $D\subseteq D'$. There are two cases.

(a) If $D'=D$, then $|D'|=|D|=d.$ By the definition of $(2\alpha'-d+1)$-disintegration, we have
\begin{align*}
N(K_\mathcal{R},Q)
&\leq \sum_{i=1}^{n-d} \sum_{j=1}^{s}\frac{P(r-r_j;\mathcal{R}\backslash \{r_j\})}{\Pi_{l=1}^{q}m_l!}\binom{2\alpha'-d+1}{r-r_j}
 \binom{n-r+r_j-i}{r_j-1}+\frac{P(r;\mathcal{R})}{\Pi_{l=1}^{q}m_l!}\binom{d}{r}\\
&= g_\mathcal{R}(n,2\alpha',2\alpha'-d+1)\\
&\leq {\rm max}\{ g_{\mathcal{R}}(n,2\alpha',k),g_{\mathcal{R}}(n,2\alpha',\alpha')\}.
\end{align*}
The third inequality follows from the condition $k\leq 2\alpha'-d+1 \leq \alpha'$ and that the function $g_\mathcal{R}(n,c,x)$ is convex for $x\in[k,\alpha'].$

(b) Otherwise $D\neq D'.$ Let $u\in V(D')\backslash V(D).$ Since $d\ge \alpha'+2.$ $u$ is not adjacent to at least one vertex $v$ in $D.$  We have $d_Q(u)+d_Q(v)\ge 2\alpha'-d+2+d-1\ge 2\alpha'+1.$ Since every nonadjacent vertices have degree sum at most $2\alpha'.$ We obtain a contradiction and the theorem is proved.
\hfill $\Box$

\section{Some corollaries}
By Theorem 1, Theorem 2, Theorem 3 and since $g_{r_1,\ldots,r_s}(n,c,k)$ is a convex function of $k$, we have the following three corollaries.

{\bf Corollary 8.} {\it Let $G$ be a 2-connected nonhamiltonian graph of order $n$ and circumference $c.$  Then
$N(K_{r_1,\ldots,r_s},G)\leq {\rm max}\{g_{r_1,\ldots,r_s}(n,c,2),\, g_{r_1,\ldots,r_s}(n,c,\lfloor c/2\rfloor)\}.$}

{\bf Corollary 9.} {\it Let $G$ be a connected nontraceable graph of order $n$ and detour order $p.$ Then
$N(K_{r_1,\ldots,r_s},G)\leq {\rm max}\{g_{r_1,\ldots,r_s}(n,p-1,1),\, g_{r_1,\ldots,r_s}(n,p-1,\lfloor (p-1)/2\rfloor)\}.$}

{\bf Corollary 10.} {\it Let $G$ be a graph on $n\geq 2\alpha'+2$ vertices with matching number $\alpha'.$ Then
$N(K_{r_1,\ldots,r_s},G)\leq {\rm max}\{g_{r_1,\ldots,r_s}(n,2\alpha',0),\, g_{r_1,\ldots,r_s}(n,2\alpha',\alpha')\}.$}

 Next we determine the maximum number of copies of $K_{r_1,\ldots,r_s}$ in $n$-vertex graphs with a given circumference or detour order. In order to obtain
  these results, we need several lemmas.

{\bf Lemma 11.} {\it For graph $G=K_t\vee(K_2+\overline{K_2})$, label the vertices in $K_t,$ $K_2,$ $\overline{K_2}$ as $\{v_1,\ldots,v_t\},$ $\{y,x\},$ $\{b,a\},$ successively. If the complete multipartite graph $K_{\mathcal{R}}$ is not a complete graph,  then the number of $K_{\mathcal{R}}$ containing $x$ in $G-\{a,b\}$ is no more than the number of $K_{\mathcal{R}}$ containing $a$ in $G.$ }

{\it Proof.} Let  $R_1,R_2,\ldots,R_s$ be the partite sets of $K_{\mathcal{R}}$ with $|R_i|=r_i$ and let $Q=G-\{a,b\}.$ We will prove that if there is a $K_{\mathcal{R}}$ containing $x$ in $Q$ then there is a $K_{\mathcal{R}}'= K_{\mathcal{R}}$ containing $a$ in $G.$ If $y \notin V(K_{\mathcal{R}})$ or $\{y, x\}$ belongs to the same class of $K_{\mathcal{R}},$ we can get $K_{\mathcal{R}}'$ by replacing $x$ to $a.$ Without loss of generality, suppose $x\in R_1$ and $y\in R_2.$ Let $R_1=\{x,v_{12},\ldots,v_{1r_1} \},$ $R_2=\{y,v_{22},\ldots,v_{2r_2} \},$ $R_j=\{v_{j1},v_{j2},\ldots,v_{jr_1} \},j\in \{3,\ldots,s\}.$ Since $K_{\mathcal{R}}$ is not a complete graph, there exists $i\in [s]$ such that $r_i\ge 2.$ We distinguish three cases.

{\bf Case 1.} $r_1\ge2,$ Let $R^1_1=\{x,a,\ldots,v_{1r_1} \},$ $R^1_2=\{v_{12},v_{22},\ldots,v_{2r_2} \},$ $R^1_j=R_j,$ $j\in\{3,\dots,s\}.$ Then the complete multipartite graph $K_{\mathcal{R}}'$ in $G$ with class $R^1_i,$ $i\in [s]$ is what we need.

{\bf Case 2.} $r_2\ge 2$ and $ r_1=1.$ Let $R^2_1=\{v_{22}\},$ $R^2_2=\{y,a,\ldots,v_{2r_2} \},$ $R^1_j=R_j,\ j\in\{3,\dots,s\}$. Then the complete multipartite graph $K_{\mathcal{R}}'$ in $G$ with class $R^2_i,i\in [s]$ is what we need.

{\bf Case 3.} $r_j\ge 3,\ j\ge3$ and $r_1=r_2=1,$ with loss of generality, suppose $r_3\ge2.$ Let $R^3_1=\{v_{31}\},$ $R^3_2=\{v_{32}\},$ $R^3_3=\{a,b,v_{33},\ldots,v_{3r_3}\},$  $R^3_j=R_j,\ j\in\{4,\dots,s\}.$ Then the complete multipartite graph $K_{\mathcal{R}}'$ in $G$ with class $R^3_i,i\in [s]$ is what we need.

Since for distinct $K_{\mathcal{R}}$, $K_{\mathcal{R}}'$ is also distinct. The proof is complete. \hfill $\Box$

{\bf Lemma 12.} {\it Let $a,b$ and $c$ be positive integers with $a\le b\le c.$ Then }
\begin{equation*}
N(K_{\mathcal{R}},K_1\vee (K_{a-1}+K_{b-1}))\le
\begin{cases}
N(K_{\mathcal{R}},K_{a+b-1}),\ \ a+b\le c+1,\\
N(K_{\mathcal{R}},K_1\vee (K_{c-1}+K_{a+b-c-1})),\ \ a+b\ge c+2.
\end{cases}
\end{equation*}
{\it Proof.} This lemma is easy to prove, so we omit the details.

Using Lemma 11, we have the following result.

{\bf Lemma 13.} {\it Let $t=\lfloor c/2\rfloor.$ Let $G$ be an $n$-vertex graph with two blocks $B_1=G(n_1,c,t)$, $B_2=G(n-n_1+1,c,t).$ Suppose the cut vertex is a dominating vertex. Then we have $N(K_{\mathcal{R}},G)\le N(K_{\mathcal{R}},G(n,c,t)).$}

Let $n-1= \alpha(c-1)+\beta$, $0\le \beta\le c-2$ and let $F(n,c)=K_1\vee(\alpha K_{c-1}+K_{\beta}).$ Define $f_{r_1,\ldots,r_s}(n,c)$ to be the number of copies of $K_{r_1,\ldots,r_s}$ in $F(n,c).$

{\bf Corollary 14.} {\it Let $G$ be an $n$-vertex nonhamiltonian graph with circumference $c.$ Then $N(K_{r_1,\ldots,r_s},G)\leq {\rm max}\{f_{r_1,\ldots,r_s}(n,c),\, g_{r_1,\ldots,r_s}(n,c,\lfloor c/2\rfloor)\}.$}

{\it Proof.} Let $n-1\ge c\ge 3,$ $t=\lfloor c/2\rfloor.$ Suppose $G$ is a graph with order $n$ and circumference $c$. If the connectivity of $G$ is at least $2,$ i.e. $G$ is $2$-connected, by Corollary 8, we have
\begin{align*}
  N(K_{\mathcal{R}},G)& \le {\rm max}\{g_\mathcal{R}(n,c,2),g_\mathcal{R}(n,c,t)\}\\
  & \le {\rm max}\{g_\mathcal{R}(n,c,1),g_\mathcal{R}(n,c,t)\}\le {\rm max}\{f_\mathcal{R}(n,c),g_\mathcal{R}(n,c,t)\}.
\end{align*}

If the connectivity of $G$ is $1,$ $G$ has at least two blocks. If $G$ has two end blocks $B_1,B_2$ such that $V(B_1)\cap V(B_2)=\emptyset,$ let $v_1$ be the cut vertex in $B_1$  and $v_2$ be the cut vertex in $B_2.$ Denote by $N_{B_1}(v_1)$ the neighbours of $v_1$ in $B_1$ and $N_{B_2}(v_2)$ the neighbours of $v_2$ in $B_2.$ Without loss of generality, suppose $d_G(v_2)-|N_{B_2}(v_2)|\ge d_G(v_1)-|N_{B_1}(v_1)|.$ Deleting edges between $N_{B_1}(v_1)$ and $v_1$ and connecting $v_2$ to $N_{B_1}(v_1),$ we obtain a graph $G'$ with $c(G)=c(G')$ and $N(K_{\mathcal{R}},G)\le N(K_{\mathcal{R}},G').$ We keep running this process until we finally get a graph, denoted by $G'',$ such that all blocks share a common vertex $v_0$. Obvious $c(G)=c(G'')$ and $N(K_{\mathcal{R}},G)\le N(K_{\mathcal{R}},G'').$ Let $B_i$, $i\in [p]$ be the blocks of $G$ with $n_i$ vertices. We have $N(K_\mathcal{R},B_i)\le {\rm max}\{f_\mathcal{R}(n_i,c),g_\mathcal{R}(n_i,c,t)\}$ by corollary 8. Since both $F(n_i,c)$ and $G(n_i,c,t)$ have dominating vertex, we suppose $v_0$ is a dominating vertex of $G$.

If for all $i\in [p]$, we have $N(K_\mathcal{R},B_i)\le f_\mathcal{R}(n_i,c),$ and hence $N(K_\mathcal{R},G)\le f_\mathcal{R}(n,c)$ by lemma 12.

If for all $i\in [p]$, we have $N(K_\mathcal{R},B_i)\le g_\mathcal{R}(n_i,c,t),$ and hence $N(K_\mathcal{R},G)\le g_\mathcal{R}(n,c,t)$ by lemma 13.

Otherwise, there exists positive integer $q$ such that $N(K_\mathcal{R},B_i)\le g_\mathcal{R}(n_i,c,t)$ for $i\in [q]$ and $N(K_\mathcal{R},B_i)\le f_\mathcal{R}(n_i,c)$ for $i\in [p]\setminus [q].$ By lemma 13, we may assume $q=1.$

If $g_\mathcal{R}(n_1,c,t)\le f_\mathcal{R}(n_1,c),$ then by lemma 12,
$$
  N(K_\mathcal{R},G)\le g_\mathcal{R}(n_1,c,t)+ f_\mathcal{R}(n-n_1+1,c)\le f_\mathcal{R}(n_1,c)+f_\mathcal{R}(n-n_1+1,c)\le f_\mathcal{R}(n,c).
$$
So we can suppose that
\begin{equation}
 g_\mathcal{R}(n_1,c,t)> f_\mathcal{R}(n_1,c).
\end{equation}
If $K_\mathcal{R}=K_s,$ one can check that $g_\mathcal{R}(n_1,c,t)\le f_\mathcal{R}(n_1,c).$ So in this case, $K_\mathcal{R}\ne K_s.$

{\bf Claim.} {\it Adding additional $c-1$ vertices to $\overline{K_{n_1-c-1+t}}$ of $G(n_1,c,t)$ will produce more number of copies of $K_{\mathcal{R}}$ than the number of copies of $K_{\mathcal{R}}$ in $K_c.$}

{\it Proof.} If $K_\mathcal{R}=S_r,r\ge3,$ this claim can be easily checked. So we assume that $K_\mathcal{R}\ne K_s$ or $K_\mathcal{R}\ne S_r.$ Let $z_1=n_1+c-1$ and $z_1-1=\alpha_1(c-1)+\beta_1$ where $2\le \alpha_1,$$ 0\le\beta_1\le c-2.$ Let $Q=G(z_1,c,t)=K_t\vee(K_{c+1-2t}+\overline{K_{z_1-c-1+t}}).$ We denote by $A=K_t,B=K_{c+1-2t},C=\overline{K_{z_1-c-1+t}}$ and order the vertices of $Q$ in $A,B,C$ with $v_1,\ldots,v_t,$ $v_{t+1},\ldots,v_{c+1-t},$ $v_{c+2-t},\ldots,v_{z_1}$ successively. Let $T_0$ be the set of the first $\beta_1+1$ vertices in $Q$ and $Q_0$ be the subgraph of $Q$ induced by $T_0.$ Then $Q_0\subseteq K_{\beta_1+1}$ and
\begin{equation}
  N(K_\mathcal{R},Q_0)\le N(K_\mathcal{R},K_{\beta_1+1}).
\end{equation}
We divide the remaining vertices of $Q$ into $\alpha_1$ set of size $c-1$ as its order. Let $T_i$ be the set of $i$-th $c-1$ vertices of $Q-Q_0$ and $Q_i$ be the subgraph induced by $\bigcup_{j=0}^iT_j.$ Let $N_i=N(K_\mathcal{R},Q_i)-N(K_\mathcal{R},Q_{i-1})$ for $i=1,2,\ldots,\alpha_1.$ Let $x_i$ be the number of $K_\mathcal{R}$ in $Q-\{v_{i+1},\ldots,v_{z_1}\}$ containing $v_i$ for $i=1,\ldots,z_1.$ It can be checked that if $c$ is even, then $t= c/2$, $B=K_1$ and $x_i$ is nondecreasing. Hence if $\alpha_1\ge 2,$ then
\begin{equation}
N_{i+1}\ge N_i,i\in [\alpha_1-1].
\end{equation}
If $c$ is odd, then $t=(c-1)/2$ and $B=K_2.$ It can be checked that $x_{i+1}\ge x_i$ for $i\in [n_1-1]\backslash \{t+2\}.$ By lemma 11 and $K_\mathcal{R}\ne K_r,$ we deduce that $x_{t+4}\ge x_{t+2}.$ So if $\alpha_1\ge 2,$ we also have $N_{i+1}\ge N_i,i\in [\alpha_1-1].$ By (1) and (2), we deduce that
$$
\sum_{i=1}^{\alpha_1-1}N_i\ge N(K_\mathcal{R},(\alpha_1-1)K_c)=(\alpha_1-1)N(K_\mathcal{R},K_c).
$$
By (3), we have $N_{\alpha_1}\ge N_{\alpha_1-1}\ge N(K_\mathcal{R},K_c).$ This completes the proof of Claim.

Let $n-n_1=\alpha_2(c-1)+\beta_2.$ Hence
$$
g_\mathcal{R}(n,c,t)\le g_\mathcal{R}(n_1,c,t)+f_\mathcal{R}(n-n_1+1,c)\le g_\mathcal{R}(n-\beta_2,c,t)+N(K_{\mathcal{R}},K_{\beta_2+1}).
$$
Recall that $K_\mathcal{R}=K_{r_1,\dots,r_s}$ and $r=r_1+\ldots+r_s.$ If $\beta_2+1<r,$ then we are done. so suppose that $r\le \beta_2+1<c.$

Note that
$$ N(K_{\mathcal{R}},K_{\beta_2+1})\le\frac{\beta_2}{c-1}N(K_{\mathcal{R}},K_c) \le \frac{\beta_2}{c-1}N_{\alpha_1}$$
and
$$ N(K_{\mathcal{R}},H(n,c,t))-N(K_{\mathcal{R}},H(n-\beta_2,c,t))\ge\frac{\beta_2}{c-1}N_{\alpha_1}.$$

Hence $$g_\mathcal{R}(n-\beta_2,c,t)+N(K_{\mathcal{R}},K_{\beta_2+1})\le g_\mathcal{R}(n,c,t).$$ If $G$ is disconnected, simply apply induction to each component of $G$ to obtain the desired result. This completes the proof. \hfill $\Box$

Next we will determine the maximum number of $K_{r_1,\ldots,r_s}$ in $n$-vertex graphs with detour order $p.$ In the following, we introduce some lemmas that will be used in our forthcoming proofs.

{\bf Lemma 15.} {\it Let $a,b$ and $p$ be positive integers with $a\le b\le p-1.$ Then }
\begin{equation*}
N(K_{\mathcal{R}},K_a)+N(K_{\mathcal{R}},K_b)\le
\begin{cases}
N(K_{\mathcal{R}},K_{a+b}),\ \ a+b\le p,\\
N(K_{\mathcal{R}},K_{p-1})+N(K_{\mathcal{R}},K_{a+b-p+1}),\ \ a+b\ge p+1.
\end{cases}
\end{equation*}

{\bf Lemma 16.} {\it Let $t=\lfloor (p-1)/2\rfloor.$ For two positive integers $n_1,n_2,$ we have
$$N(K_{\mathcal{R}},G(n_1,p-1,t))+N(K_{\mathcal{R}},G(n_2,p-1,t))\le N(K_{\mathcal{R}},G(n_1+n_2,p-1,t)).$$}
Let $n= \alpha p+\beta$, $0\le \beta\le p-1$. Define $h_{r_1,\ldots,r_s}(n,p)$ be the number of $K_{r_1,\ldots,r_s}$ in $\alpha K_p +K_{\beta}$ i.e. $h_{r_1,\ldots,r_s}(n,p)=N(K_{r_1,\ldots,r_s},\alpha K_p +K_{\beta}).$

{\bf Corollary 17.} {\it Let $G$ be a nontraceable graph of order $n$ and detour order $p.$ Then
$$
N(K_{r_1,\ldots,r_s},G)\leq {\rm max}\{h_{r_1,\ldots,r_s}(n,p),\, g_{r_1,\ldots,r_s}(n,p-1,\lfloor (p-1)/2\rfloor)\}.
$$}
{\it Proof.} Let $n-1\ge p\ge 3,$ $t=\lfloor (p-1)/2\rfloor$ and let $G$ be a graph with order $n$ and detour order $p$. If $G$ is connected, by Corollary 12, we have
\begin{align*}
  N(K_{\mathcal{R}},G)\le& {\rm max}\{g_\mathcal{R}(n,p-1,1),g_\mathcal{R}(n,p-1,t)\}\\
  \le& {\rm max}\{g_\mathcal{R}(n,p-1,0),g_\mathcal{R}(n,p-1,t)\}\le {\rm max}\{h_\mathcal{R}(n,p),g_\mathcal{R}(n,p-1,t)\}.
\end{align*}

 Otherwise, suppose $G$ is not connected. Let $G_i$, $i\in [l],$ be the components of $G$ with $n_i$ vertices. We have $N(K_\mathcal{R},G_i)\le {\rm max} \{h_\mathcal{R}(n_i,p),g_\mathcal{R}(n_i,p-1,t)\}.$

If for all $i\in [l]$, we have $N(K_\mathcal{R},G_i)\le h_\mathcal{R}(n_i,p),$ and then by Lemma 15, we deduce that $N(K_\mathcal{R},G)\le h_\mathcal{R}(n,p).$

If for all $i\in [l]$, we have $N(K_\mathcal{R},G_i)\le g_\mathcal{R}(n_i,p-1,t),$ and then by Lemma 16, we deduce that $N(K_\mathcal{R},G)\le g_\mathcal{R}(n,p-1,t).$

Otherwise, we suppose $N(K_\mathcal{R},G_i)\le g_\mathcal{R}(n_i,p-1,t)$ for $i\in [q]$ and $N(K_\mathcal{R},G_i)\le h_\mathcal{R}(n_i,p)$ for $i\in [l]\setminus [q],$ then by lemma 16, we assume $q=1.$

If $g_\mathcal{R}(n_1,p-1,t)\le h_\mathcal{R}(n_1,p),$ then by Lemma 15,
\begin{align*}
 N(K_\mathcal{R},G)&\le g_\mathcal{R}(n_1,p-1,t)+ h_\mathcal{R}(n-n_1,p)\\
 &\le h_\mathcal{R}(n_1,p)+h_\mathcal{R}(n-n_1,p)\le h_\mathcal{R}(n,p).
\end{align*}
So we can suppose that
\begin{equation}
 g_\mathcal{R}(n_1,p-1,t)> f_\mathcal{R}(n_1,p-1).
\end{equation}
If $K_\mathcal{R}=K_s,$ one can check that $g_\mathcal{R}(n_1,p-1,t)\le f_\mathcal{R}(n_1,p-1).$ So in this case, $K_\mathcal{R}\ne K_s.$ Let $z_1=n_1+p=\alpha_1p+\beta_1$ where $2\le \alpha_1, 1\le \beta_1\le p-1$ and let $Q=G(z_1,p-1,t)=K_t\vee(K_{p-2t}+\overline{K_{z_1-p+t}}).$  We denote by $A=K_t,B=K_{p-2t},C=\overline{K_{z_1-p+t}}.$ We will prove  that adding additional $p$ vertices to $C$ of $Q$ will produce more numbers of $K_{\mathcal{R}}$ than $K_p.$ Denote the vertices of $Q$ in $A,B,C$ by $v_1,\ldots,v_t,$ $v_{t+1},\ldots,v_{p-t},$ $v_{p+1-t},\ldots,v_{z_1}$ successively. Let $T_0$ be the set of the first $\beta_1$ vertices in $Q$ and $Q_0$ be the subgraph of $Q$ induced by $T_0.$ Then $Q_0\subseteq K_{\beta_1}$ and
\begin{equation}
  N(K_{\mathcal{R}},Q_0)\le N(K_{\mathcal{R}},K_{\beta_1}).
\end{equation}

We divide the vertices of $Q-T_0$ into $\alpha_1$ set of size $p$ as its order. Let $T_i$ be the set of $i$-th $p$ vertices of $Q-Q_0$ and let $Q_i$ be the subgraph induced by $\bigcup_{j=0}^iT_j.$ Let $N_i=N(K_{\mathcal{R}},Q_i)-N(K_{\mathcal{R}},Q_{i-1})$ for $i=1,2,\ldots,\alpha_1$ and $x_i$ be the number of $K_{\mathcal{R}}$ in $Q- \{v_{i+1},\ldots,v_{z_1}\}$ containing $v_i,$ for $i=1,\ldots,z_1.$ It can be checked that if $p$ is odd, then $t= (p-1)/2$, $B=K_1$ and $x_i$ is nondecreasing. Hence if $\alpha_1\ge 2,$ then
\begin{equation}
N_{i+1}\ge N_i,i\in [\alpha_1-1].
\end{equation}
If $p$ is even, then $t=(p-2)/2$ and $B=K_2.$ It can be checked that $x_{i+1}\ge x_i$ for $i\in [z_1-1]\backslash \{t+2\}.$ By Lemma 11, we deduce that $x_{t+4}\ge x_{t+2}.$ So if $\alpha_1\ge 2,$ we also have $N_{i+1}\ge N_i,i\in [\alpha_1-1].$ By (4) and (5), we have that
\begin{equation}
\sum_{i=1}^{\alpha_1-1}N_i\ge N(K_{\mathcal{R}},(\alpha_1-1)K_p)=(\alpha_1-1)N(K_{\mathcal{R}},K_p).
\end{equation}
By (6), we have $N_{\alpha_1}\ge N_{\alpha_1-1}\ge N(K_{\mathcal{R}},K_p).$ Furthermore, adding additional $p$ vertices to $C$ of $Q$ will produce more copies of $K_{\mathcal{R}}$ than $K_p.$ Let $n-n_1=\alpha_2p+\beta_2.$ Hence
$$
g_\mathcal{R}(n_1,p-1,t)+h_\mathcal{R}(n-n_1,p)\le g_\mathcal{R}(n-\beta_2,p-1,t)+N(K_{\mathcal{R}},K_{\beta_2}).
$$
By the same argument as in the proof of of Corollary 14, we have
$$g_\mathcal{R}(n-\beta_2,p-1,t)+N(K_{\mathcal{R}},K_{\beta_2})\le g_\mathcal{R}(n,p-1,t).$$
This completes the proof.\hfill $\Box$

\section{Concluding remarks}

In \cite{11}, Erd\H{o}s and Gallai determined the maximum number of edges in $n$-vertex graph with circumference $c.$ They also showed that the extremal graphs are $k_1\vee (\frac{n-1}{k-1} K_{c-1}),$ where $k-1$ divides $n-1$. In \cite{12}, Luo showed that these graphs are also the extremal examples for the maximum number of cliques in $n$-vertex graphs with circumference $c$. But as we show in Corollary 13, if $n$ is sufficiently large, the extremal graphs for the maximum number of complete multipartite graph with diameter $2$ in $n$-vertex graphs with given circumference are not the same as in the above two cases.

{\bf Acknowledgement.} The author would like to thank Prof. Xingzhi Zhan for his constant support and guidance and Yuxuan Liu for conducive discussions and careful reading of a draft. This research  was supported by the NSFC grant 11671148 and Science and Technology Commission of Shanghai Municipality (STCSM) grant 18dz2271000.

\end{document}